\documentclass[final,reqno,11pt,centertags]{article}
\usepackage{amsmath,amsthm,amscd,amssymb,latexsym,upref}
\date{\today}
\input epsf
\usepackage{epsfig}
\usepackage[T2A,OT1]{fontenc}
\usepackage[ot2enc]{inputenc}
\usepackage[russian,english]{babel}

\usepackage{epic,eepicemu}

\newcommand{\bbD}{{\mathbb{D}}}

\newcommand{\bbR}{{\mathbb{R}}}

\newcommand{\bbZ}{{\mathbb{Z}}}
\newcommand{\bbC}{{\mathbb{C}}}

\newcommand{\bbT}{{\mathbb{T}}}

\newcommand{\cM}{{\mathcal{M}}}

\newcommand{\fA}{{\mathfrak{A}}}

\renewcommand{\Re}{\text{\rm Re}\,}
\newcommand{\Res}{\text{\rm Res}}
\renewcommand{\Im}{\text{\rm Im}\,}

\newcommand{\tr}{\text{\rm tr}}

\allowdisplaybreaks \numberwithin{equation}{section}
\newtheorem{theorem}{Theorem}[section]

\newtheorem{lemma}[theorem]{Lemma}
\newtheorem{proposition}[theorem]{Proposition}
\newtheorem{corollary}[theorem]{Corollary}

\theoremstyle{definition}
\newtheorem{definition}[theorem]{Definition}

\newtheorem{remark}[theorem]{Remark}
\newtheorem{problem}[theorem]{Problem}

\newtheorem{example}[theorem]{Example}

\title{On the Direct Cauchy Theorem in Widom Domains: Positive and Negative Examples}
\author{ P. Yuditskii\thanks{
Supported by the Austrian Science Fund FWF, project no: P22025--N18}}

\date{\today}
\begin{document}
\maketitle

%\textit
\begin{abstract}
We discuss several questions which remained open in our joint work with M. Sodin
"Almost periodic Jacobi matrices
with homogeneous spectrum, infinite-dimensional Jacobi inversion,
and Hardy spaces of character--automorphic functions". In particular, we show that there exists a non-homogeneous set $E$ such that  the Direct Cauchy Theorem (DCT) holds in the Widom  domain $\bbC\setminus E$.
On the other hand we demonstrate that the weak homogeneity condition on $E$ (introduced recently by
Poltoratski and  Remling)  does not ensure that DCT holds
in the corresponding Widom  domain.

\end{abstract}

\section{Introduction}
Several recent publications \cite{AM, BRZ, CSZ, GZ, PR, SIM} indicate a certain interest in our paper \cite{SYU}. We recall the main result of this paper, simultaneously we introduce notations and give definitions necessary in what follows.

Let $E$ be a compact on the real axis without isolated points, $E=[b_0,a_0]\setminus \cup_{j \ge 1}(a_j,b_j)$.
By $J(E)$ we denote the set of reflectionless (two-sided) Jacobi matrices $J$ with the spectrum on $E$. This means that the diagonal elements of the resolvent
\begin{equation*}
    R_{k,k}(z)=\langle(J-z)^{-1}e_k,e_k\rangle
\end{equation*}
possess the property $\Re R_{k,k}(x+i0)=0$ for almost all $x\in E$. As usual $e_k$'s denote the standard basis in $l^2$.
For an exceptional role of this class of Jacobi matrices see \cite{CR}.

The function $R(z)=R_{k,k}(z)$ has positive imaginary part in the upper half plane, and therefore possesses the representation
\begin{equation}\label{intreprR}
    R(z)=\int_E\frac{d\sigma(x)}{x-z}.
\end{equation}
Moreover, since $R(z)$ assumes \textit{pure imaginary values} on $E$ we have
$$
\sigma'_{a.c.}=\frac 1\pi|R(x)|.
$$
We follow the terminology in \cite{PR} and call reflectionless the measures $\sigma$ related
to reflectionless  functions $R(z)$ \eqref{intreprR},  $\Re R(x+i0)=0$, a.e. $x\in E$.

The collection of reflectionless functions associated to the given compact $E$ can be parameterized in the following way.    We chose arbitrary
$x_j\in [a_j,b_j]$ and set
\begin{equation}\label{defchi}
   R(z)= R(z,\{x_j\})=-\frac{1}{\sqrt{(z-a_0)(z-b_0)}}\prod_{j\ge 1}\frac{z-x_j}{\sqrt{(z-a_j)(z-b_j)}}.
\end{equation}

By $D(E)$ we denote the set of so-called divisors $D$, where
\begin{equation}\label{divisors}
    D=\{(x_j,\epsilon_j): \ x_j\in[a_j,b_j],\ \epsilon_j=\pm 1\},\quad (a_j,1)\equiv(a_j,-1),\
   (b_j,1)\equiv(b_j,-1).
\end{equation}
The map $J(E)\to D(E)$ is defined in the following way. For a reflectionless $J$ the resolvent function $R_{0,0}(z)$ possesses the representation \eqref{defchi} and this representation produces the collection $\{x_j\}$. To define $\epsilon_j$ we represent $J$ as a two-dimensional perturbation of the block-diagonal sum of one-sided Jacobi matrices $J_\pm$, that is,
\begin{equation*}
    J=\begin{bmatrix}
    J_-&0\\
    0&J_+
    \end{bmatrix}+p_0e_{-1}\langle\cdot, e_0\rangle+p_0e_{0}\langle\cdot, e_{-1}\rangle.
\end{equation*}
This representation generates the identity
\begin{equation}\label{sumrs}
    -\frac{1}{R_{0,0}(z)}=-\frac{p_0^2}{ r_-(z)}+r_+(z),
\end{equation}
where
\begin{equation*}
    r_+(z)=\langle(J_+-z)^{-1} e_{0},e_{0}\rangle,\quad r_-(z)=\langle(J_--z)^{-1} e_{-1},e_{-1}\rangle.
\end{equation*}
For $x_j\in(a_j,b_j)$ this point is a pole of only one of the two functions in the right hand side of \eqref{sumrs}. We set $\epsilon_j=1$  if $x_j$ is a pole of $r_+(z)$ and $\epsilon_j=-1$ in the opposite case.

Next we will define the so called generalized Abel map, the map from the collection of divisors $D(E)$ to the group of characters of the fundamental group of the domain $\Omega=\bar\bbC\setminus E$. We are able to do this for  Widom domains.

Let $z:\bbD/\Gamma\to\Omega$ be a uniformization of the domain $\Omega$, that is, an analytic function $F$ in $\Omega$ can be represented by an analytic function $f$ in $\bbD$, which is automorphic with respect to the action of the Fuchsian group $\Gamma$:
\begin{equation*}
    F(z(\zeta))=f(\zeta),\ f(\gamma(\zeta))=f(\zeta),\ \zeta\in\bbD,\ \gamma\in\Gamma.
\end{equation*}
The dual group $\Gamma^*$ is formed by characters
\begin{equation*}
    \alpha:\Gamma\to\bbT\ \text{such that}\ \alpha(\gamma_1\gamma_2)=\alpha(\gamma_1)\alpha(\gamma_2),
    \ \gamma_{1,2}\in\Gamma.
\end{equation*}
For a character $\alpha$ the Hardy space $H^p(\alpha)$ is the subspace of the standard $H^p$ consisting of character automorphic functions $f(\gamma(\zeta))=\alpha(\gamma)f(\zeta)$.

A domain is called of Widom type if $H^\infty(\alpha)$ is not trivial (contains a non constant function) for all $\alpha\in\Gamma^*$. Let $b(\zeta,\zeta_0)$ be the Green function of the group $\Gamma$ \cite{POM}, i.e.,
\begin{equation}\label{defcgre}
    b(\zeta,\zeta_0)=\prod_{\gamma\in\Gamma}\frac{\gamma(\zeta_0)-\zeta}{1-\zeta\overline{\gamma(\zeta_0)}}
    \frac{|\gamma(\zeta_0)|}{\gamma(\zeta_0)}
\end{equation}
Note that $G(z(\zeta),z_0):=-\log |b(\zeta,\zeta_0)|$ is the Green function in the domain $\Omega$. We assume that $z(0)=\infty$, that is, $G(z)=-\log |b(\zeta,0)|$ is the Green function with respect to infinity.

A regular domain $\Omega$ is of Widom type if and only if
\begin{equation}\label{widomcond}
    \sum_{c_k:\nabla G(c_k)=0}G(c_k)<\infty.
\end{equation}
Let us point out that  the critical point $c_k$ belongs to $(a_k,b_k)$. The condition \eqref{widomcond} guaranties that the following product
\begin{equation}\label{kade}
    K^D(\zeta)=\sqrt{\prod_{j\ge 1}\frac{z(\zeta)-x_j}{z(\zeta)-c_j}\frac{b(\zeta,c_j)}{b(\zeta,x_j)}}
    \prod_{j\ge 1}b(\zeta,x_j)^{\frac{1+\epsilon_j}{2}}
\end{equation}
converges for an arbitrary $D\in D(E)$.

The Abel map $D(E)\to\Gamma^*$ is defined by the relation
\begin{equation*}
    K^D(\gamma(\zeta))=\alpha^D(\gamma)K^D(\zeta).
\end{equation*}

The third map $\Gamma^*\to J(E)$ is defined as follows. Let $k^\alpha$ be the reproducing kernel in $H^2(\alpha)$ with respect to the origin, i.e., $\langle f, k^\alpha\rangle=f(0)$ for all $f\in H^2(\alpha)$.
Let $\alpha_0$ be the character of the Green function $b$, $b\circ\gamma=\alpha_0(\gamma)b$.
Then $J(\alpha)$ is the matrix of the multiplication operator by $z(\zeta)$ with respect to the orthonormal basis $\{e(\zeta,n)\}_{n\in\bbZ}$ in $L^2(\alpha)$, where
\begin{equation*}
    e(\zeta,n)=b^n\frac{k^{\alpha\alpha_0^{-n}}}{\sqrt{{k^{\alpha\alpha_0^{-n}}}(0)}}.
\end{equation*}

The main result in \cite{SYU} claims that under a certain additional condition on the domain $\Omega$ all three maps
\begin{equation}\label{allthree}
    J(E)\to D(E)\to\Gamma^*\to J(E)
\end{equation}
are one-to-one and continuous with respect to the operator norm topology in $J(E)$; recall that  $\Gamma^*$ is a compact Abelian group; $D(E)$ is equipped with the product topology.

This additional condition is called the \textit{Direct Cauchy Theorem} (DCT) \cite{HA}. In fact, it is not a theorem, but a certain property of a Widom domain. For some of them DCT holds true, for others it fails.

We say that $F(z)$ is of Smirnov class in $\Omega$ if the corresponding function $f(\zeta)=F(z(\zeta))$ is of Smirnov class in $\bbD$, that is, it possesses  a representation $f=\frac{f_1}{f_2}$, where $f_1, f_2$ are uniformly bounded, moreover the denominator is an outer function.

\begin{definition}
The space $E_0^1(\Omega)$ is formed by Smirnov class functions $F$ in
$\Omega$ such that
$$
\|F\|=\frac 1 {2\pi}\int_E |F(x+i0)|\,{dx}+\frac 1 {2\pi}\int_E |F(x-i0)|\,{dx}<\infty,
$$
and $F(\infty)=0$.
\end{definition}

\begin{definition}
We say that the \textit{Direct Cauchy Theorem (DCT)} holds if
\begin{equation}\label{dct}
    \frac 1{2\pi i}\oint_E F(x) dx=\Res_{\infty} F(z)dz= A,\quad
    F(z)=-\frac A z+\dots,\ z\to\infty,
\end{equation}
for all $F\in E_0^1(\Omega)$.
\end{definition}

\begin{theorem}[Sodin-Yuditskii \cite{SYU}] Let $E$ be such that $\Omega=\bar\bbC\setminus E$ is of Widom type with DCT. Then every Jacobi matrix $J\in J(E)$ is almost periodic.
\end{theorem}

The following notion was introduced by Carleson \cite{CAR}. A compact $E$ is homogeneous if there exists $\eta=\eta(E)>0$ such that
\begin{equation}\label{hocond}
    |E\cap(x-\delta,x+\delta)|\ge \eta\delta
\end{equation}
for all $x\in E$ and $\delta\in(0,1)$.

It was noted in \cite{SYU} that if $E$ is homogeneous then $\Omega$ is a Widom domain with DCT. Thus the homogeneity is a very nice constructive condition
on $E$
that guaranties that $J(E)$ consists of almost periodic operators. For instance all standard Cantor sets of positive length are homogeneous, for a proof see e.g. \cite{PYU}.

\smallskip
In this note we answer several remaining open questions in the above context.

\smallskip

First, the map $J(E)\to D(E)$ is one-to-one if and only if all reciprocal  $-1/R(z)$ to reflectionless  (Nevanlinna class) functions  have no singular component on the set $E$ in their integral representation, i.e., for
\begin{equation}\label{resipr}
    -\frac 1{R(z)}=z-q_0+\sum_{x_j\in(a_j,b_j)}\frac{\tau_j}{x_j-z}+\int_{E}\frac{d\tau(x)}{x-z}
\end{equation}
the singular component of the measure $\tau$ is trivial, $\tau_s(E)=0$.

The fact that neither reflectionless measure nor its reciprocal \eqref{resipr} has no singular component on $E$ for Widom domains with DCT was proved in \cite{PYU}.

In general, the question on the possible support of the singular part of a reflectionless measure was studied in \cite{PR}. In particular they introduced a notion of a weakly homogeneous set: a Borel set $E$ is weakly homogeneous if
\begin{equation}\label{whom}
    \limsup_{\delta\to+0}\frac 1\delta  |E\cap(x-\delta,x+\delta)|>0,\quad x\in E,
\end{equation}
and proved that all reflectionless measures on a weakly homogeneous set are absolutely continuous.

Moreover, previously known examples of Widom domains such that DCT fails were based on the idea to construct a reflectionless measure with a non-trivial singular component,   for details see Sect. \ref{sec5}.

So,  assume a priory that  $\Omega=\bar\bbC\setminus E$ is of Widom type.
\begin{itemize}
\item
Does the weak homogeneity \eqref{whom} imply DCT in this case?
\end{itemize}

Note that this assumption \eqref{whom} is even stronger than
%to try to deduce DCT from
the property that all reflectionless measures  associated with the given $E$ (and their reciprocals \eqref{resipr}) are absolutely continuous.

Second, the map $D(E)\to \Gamma^*$, which was defined for Widom domains, is one-to-one if and only if DCT holds. It deals with the following property of $H^2$-spaces in $\Omega$.

In the classical Hardy spaces theory there is an important description of their orthogonal complement in the standard $L^2$, namely, $\bar f\in L^2\ominus H^2$ if and only if $f\in H^2$ and $f(0)=0$. In the Widom domain case the following statement holds true. Define the Blaschke product
\begin{equation*}
    \theta(\zeta)=\prod_{j\ge 1}b(z,c_j),
\end{equation*}
which converges due to the Widom condition \eqref{widomcond}, and denote by $\alpha_\theta$ the corresponding character, $\theta\circ\gamma=\alpha_\theta(\gamma)\theta$. Then, if $\bar f\in L^2(\alpha)\ominus H^2(\alpha)$, then $\theta f\in H^2(\alpha_\theta\alpha^{-1})$ and $f(0)=0$.
In other words, if we define
\begin{equation*}
    \check H^2(\alpha)=\{f\in L^2(\alpha): \ \theta\bar f\in L^2(\alpha_\theta\alpha^{-1})\ominus
    H_0^2(\alpha_\theta\alpha^{-1})\},
\end{equation*}
then $\check H^2(\alpha)\subset H^2(\alpha)$.

Thus, generally speaking, for a Widom domain to the given character $\alpha$ one can associate the biggest $H^2(\alpha)$ and the smallest $\check H^2(\alpha)$ possible  Hardy spaces. They coincide, i.e., $H^2(\alpha)=\check H^2(\alpha)$ for all $\alpha\in\Gamma^*$, if and only if DCT holds in the Widom domain \cite{SYU}. Both reproducing kernels $k^{\alpha}/\|k^{\alpha}\|$ and $\check k^\alpha/\|\check k^\alpha\|$ are functions of the form \eqref{kade}. Thus, as soon
as $H^2(\alpha)\not=\check H^2(\alpha)$ two different divisors in $D(E)$ correspond to the same character $\alpha$.

\begin{itemize}
\item
How big can the defect subspace $H^2(\alpha)\ominus\check H^2(\alpha)$ be?
\end{itemize}

Finally, it is worthwhile to clarify:
\begin{itemize}
\item
Is the homogeneity \eqref{hocond} just a sufficient condition
  for DCT or  is it also  necessary?
\end{itemize}

To answer these three questions we first relate DCT with an  $L^1$-extremal problem, Sect. \ref{secl1}. In Sect. \ref{secstru} we reveal the structure of the extremal function, Theorem \ref{threpr}.  It shows, that in the simplest case, when $E$ consists of a system of intervals having a unique accumulation point $x_0\in E$,  our questions can be reduced  to approximation problems  for \textit{entire functions} (with respect to the variable $1/(z-x_0)$), e.g., see Lemma \ref{lemma 7.1}. This area was developed essentially recently by  Borichev-Sodin \cite{BoS1, BoS2, BoS3}.

We can summarize the results of the last two sections in the following proposition.

\begin{proposition}
Define the following three classes of  Widom(-Denjoy) domains:
 \begin{itemize}
 \item
 $\Omega=\bar\bbC\setminus E\in W_{hom}$ if $E$ is homogeneous;
   \item
   $\Omega\in W_{DCT}$ if Direct Cauchy Theorem holds in $\Omega$;
   \item
   $\Omega\in W_{a.c}$ if all reflectionless measures \eqref{intreprR} in $\Omega$ and their reciprocal
   \eqref{resipr} are absolutely continuous on $E$.
 \end{itemize}
Then
\begin{equation}\label{www}
    W_{hom}\subset W_{DCT}\subset W_{a.c.}
\end{equation}
and both  inclusions are   proper.
\end{proposition}

\bigskip
\noindent
\textbf{Acknowledgment.} Initial push to consider these problems was generated by a small discussion with Barry Simon. Then we discussed the subject with Sasha Volberg, the manuscript \cite{NVYu} is the result of this discussion. I have to mention a special role of Misha Sodin: I am very thankful for his explanations
dealing with the potential
theory in the Denjoy domains, in particular of Benedicks' construction
used in section \ref{secbndx}.

% The special role of Misha Sodin was mentioned above.

\section{$L^1$ extremal problem in Widom Domain \\ and Direct Cauchy Theorem}\label{secl1}

We consider here the following extremal problem.
\begin{problem}\label{pr1} Find
\begin{equation}\label{inf}
    M=\inf\{\|F\|: F\in E_0^1(\Omega),\ \ F(z)=-\frac 1 z+\dots,\ z\to\infty\}.
\end{equation}
\end{problem}
This extremal problem is closely related with the DCT.

\begin{theorem}\label{thm1}
 For a Widom domain $\Omega$ the DCT holds if and only if  $M=1$ in Problem \ref{pr1}.
\end{theorem}

\begin{proof}%[Proof of Theorem \ref{thm1}]
We adapt the general proof of this kind of theorem \cite{HA} to the special Denjoy domain case.
Let $R(z)$ be a reflectionless function, see \eqref{defchi}.
It was shown in \cite{SYU} that $R(z)$ is of Smirnov class, so $R\in E^1_0(\Omega)$ and we get
\begin{equation}\label{pone}
    M\le \|R\|\le \int_E d\sigma=1.
\end{equation}
On the other hand, if DCT holds, then for every function $F$ of the set \eqref{inf} we have
\begin{equation*}
    1=\frac 1{2\pi i}\oint_E F\le \|F\|.
\end{equation*}
Thus $M=1$.

Conversely, $M=1$  implies
$$
|\Lambda(F)|\le \|F\|
$$
for the functional
$$
\Lambda(F)=A,\ F\in E^1_0(\Omega), \ F=-\frac A z+\dots,\ z\to\infty.
$$

Since $E^1\subset L^1$ we can extend $\Lambda$ to a functional in $L^1$. Therefore there exists $w\in L^\infty$,
$\|w\|\le 1$ such that
\begin{equation}\label{funcA}
    A=\frac 1{2\pi i}\oint_E w(z)F(z)dz
\end{equation}

Now we put here $F=R$. We get
\begin{equation*}
\begin{split}
1=&\frac 1{2\pi i}\oint_E w(z)R(z)dz\\
=& \frac 1{2\pi }\oint_E w(z)|R(z)dz|
\le\frac 1{2\pi }\oint_E |w(z)|\,|R(z)dz|
\\ \le&\frac 1{2\pi }\oint_E\,|R(z)dz|\le 1
\end{split}
\end{equation*}
Thus $w=1$, and hence, by \eqref{funcA},
$$
A=\frac 1{2\pi i}\oint_E F(z)dz
$$
for all $F\in E^1$, $F(\infty)=0$.

\end{proof}

%$$
%\oL^2_{|R|}=L^2_{|R|}\oplus R^{-1}L^2_{|R|^{-1}}
%$$

In the end of this section we show that the infimum \eqref{inf} is assumed.

Except for \eqref{widomcond} one can characterize a Widom domain by the following property.
Let $\omega(dx,z)$ be the harmonic measure in the domain with respect to  $z\in\Omega$.
Then
the harmonic measure
 $\omega(dx):=\omega(dx,\infty)$ is absolutely continuous $\omega(dx)=\rho(x)\, dx$, moreover \cite{POM}
\begin{equation*}
    \int _E\rho(x)\log\rho(x)\,dx>-\infty.
\end{equation*}
So we can define an outer  function $\Phi(z)$ multi--valued in $\Omega$, such that
\begin{equation*}
    \frac1 {{2\pi}}|{\Phi(z)}|=e^{\int_E\log\rho(x)\omega(dx,z)},
\end{equation*}
i.e., $\rho(x)= \frac 1 {2\pi}|\Phi(x)|$, $|\Phi(x)|=\lim_{\epsilon\to 0}|\Phi(x\pm \epsilon i)|$ for a.e. $x\in E$.
By $\alpha_1$ we denote the character generated by $\Phi^{-1}(z(\zeta))$,
$$
\Phi^{-1}(z(\gamma(\zeta)))=\alpha_1(\gamma)\Phi^{-1}(z(\zeta)).
$$

Let us use an alternative description of
the space $E^1_0(\Omega)$.

\begin{proposition}\label{pr2.3} Let
$$
 H_0^1(\alpha_1):=\{f\in H^1: f(\gamma(\zeta))=\alpha_1(\gamma)f(\zeta),\ f(0)=0\}.
$$
$f\in  H_0^1(\alpha_1)$ if and only if $f(\zeta)=F(z(\zeta))$ and
$F\Phi\in E^1_0(\Omega)$.
\end{proposition}

\begin{proof}
It follows from the fact that one can identify $f\in H^1(\alpha_1)$ with a multi--valued
Smirnov class function $F(z)$, $F(z(\zeta))=f(\zeta)$, which is integrable with respect to the
harmonic measure
$$
\int_{\bbT}|f(t)|dm(t)=\oint_E |F(x)|\omega(dx)=\frac 1 {2\pi}\oint_E |F(x)\Phi(x)||dx|.
$$
Note also that $F\Phi$ is singlevalued in $\Omega$ by the definition of the character $\alpha_1$.
\end{proof}

\begin{lemma} There exists $H\in E^1_0(\Omega)$ such that
$$
M=\|H\|.
$$
\end{lemma}
\begin{proof}
In an extremal sequence we chose a subsequence that converges
uniformly on compact subsets in $\Omega$, that is,
$$
H(z)=\lim F_n(z), \quad z\in \Omega.
$$
In other words $f_n(\zeta)=F_n\Phi^{-1}(z(\zeta))$ converges to $h(\zeta)=H\Phi^{-1}(z(\zeta))$ uniformly on  compact subsets in $\bbD$.
Therefore
$$
\int_{\bbT} |h(rt)|dm(t)=\lim_{n\to\infty}\int_{\bbT} |f_n(rt)|dm(t)\le \lim_{n\to\infty}\|f_n\|.
$$
Thus $h\in H^1_0(\alpha_1)$, moreover $\|h\|\le \lim_{n\to\infty}\|f_n\|$. In other words
$H\in E^1_0(\Omega)$ and $\|H\|\le \lim_{n\to\infty}\|F_n\|$.

Therefore we get
\begin{equation*}
\inf\|F(x)\|\le
\|H(x)\|\le \lim_{n\to\infty}\|F_n\|=
\inf\|F(x)\|.
\end{equation*}

\end{proof}

\section{Reduction to a weighted $L^2$ extremal function. Structural Theorem}\label{secstru}
\begin{lemma} Among the extremal functions of Problem \ref{pr1} there is a function of the form
\begin{equation}\label{exfunc2}
H(z)=-\frac{1}{\sqrt{(z-a_0)(z-b_0)}}\prod_{j\ge 1}\sqrt{\left(\frac{z-a_j}{z-b_j}\right)^{\tilde\delta_j}} K^2(z), \ \tilde\delta_j=\pm 1,
\end{equation}
where $K$ is a single--valued Smirnov class function, $K(\infty)=1$.
\end{lemma}

\begin{proof}
We can assume that $H(z)\in \bbR$,
for $z\in \bbR\setminus E$, otherwise we
use the extremal function $\frac{H(z)+\overline{H(\bar z)}}2$.
Let us prove some properties of such a function $H$.

$H(z)$ has not more than one zero in each gap $(a_j,b_j)$.
Indeed, if $H(z_1)=H(z_2)=0$, $a_j<z_1<z_2<b_j$, then the function
$$
H(z)\left(1-\frac{\epsilon}{(z-z_1)(z-z_2)}\right)
$$
has the smaller $E^1$-norm for a sufficiently small $\epsilon$.

In the same way we can show that the extremal function has no
complex zeros and also $H(z)\not= 0$ for $z\in (a_0,\infty)\cup(-\infty,b_0)$.

Thus the extremal function is of the form (see \eqref{defchi})
\begin{equation}\label{exfunc1}
H(z)=R(z,\{x_j\}) F(z),
\ F(\infty)=1,
\end{equation}
where $x_j\in[a_j,b_j]$ and $F(z)$ is of Smirnov class, $F(z)\not=0$ in $\bbC\setminus E$ and such that
\begin{equation}\label{intxj}
\begin{split}
   \frac 1 {2\pi}& \oint_{E}|F(x)||R(x,\{x_j\})|\, dx\\=
   \frac 1 {\pi}&\int_E |F(x)|\prod_{j\ge 1}\frac{x-x_j}{\sqrt{(x-a_j)(x-b_j)}}\frac{dx}{\sqrt{(a_0-x)(x-b_0)}}<\infty.
    \end{split}
\end{equation}

Assume now, that  $x_k\in(a_k,b_k)$ for a certain $k$, i.e., it is not one of the end points. For all other parameters frizzed (including the function $F$)
we consider the integral in \eqref{intxj} as a function of the given $x_k$. It is well defined in this domain ($x_k\in(a_k,b_k)$),
 moreover it represents  a linear function  of $x_k$. Therefore the infimum is assumed  on the left or  right boundary point  of the interval.
 That
is, the extremal  function is of the form
 \begin{equation}\label{exfunc3}
H(z)=-\frac{1}{\sqrt{(z-a_0)(z-b_0)}}\prod_{j\ge 1}\sqrt{\left(\frac{z-a_j}{z-b_j}\right)^{\delta_j}} F(z), \ \delta_j=\pm 1.
\end{equation}

The function $F$ has no zeros and it is single--valued in the domain. Let $\gamma_j$ be a contour that starts at the "upper" bound of the interval
$(a_j,b_j)$, and goes through infinity to the "lower" bound.
Then the change of its argument along the contour is of the form
$$
\Delta_{\gamma_j}\arg F=2\pi n_j.
$$
Now we represent $F$ in the form
$$
 F(z)=\prod_{\{j: n_j\ {\rm is\ odd} \}}{\left(\frac{z-a_j}{z-b_j}\right)^{-\delta_j}}\tilde F(z)
 $$
 The point is that $\sqrt {\tilde F(z)}$ is single--valued in $\Omega$ and we set $K(z)=\sqrt {\tilde F(z)}$.
Thus the lemma is  proved.
\end{proof}

The following lemma is almost evident. For a given weight $w$ we define $E^2_w= E^2_w(\Omega)$ as a set of single--valued Smirnov class functions,
which are square--integrable against $w|dx|$, i.e.:
\begin{equation}\label{defhw}
\|F\|^2_w=\frac 1{2\pi}\oint _E|F|^2 w|dx|<\infty.
\end{equation}

\begin{lemma} \label{l7} Let $\chi(x)$ be the weight function of the form
\begin{equation}\label{defwchi}
   \chi(x)= \frac{1}{\sqrt{(a_0-x)(x-b_0)}}\prod_{j\ge 1}\sqrt{\left(\frac{x-a_j}{x-b_j}\right)^{\tilde\delta_j}},
\end{equation}
which corresponds to the particular choice of $x_j$ given by
\eqref{exfunc2}.
Let $k$ be the reproducing kernel in this space with respect to $\infty$.  Then
\begin{equation}\label{hask}
    K(z)=\frac{k(z)}{k(\infty)}.
\end{equation}

\end{lemma}

\begin{proof}
For every function $F\in E^2_{\chi}$, $F(\infty)=0$, we have
$$
\min_{\epsilon}\oint_E|K(x)+\epsilon F(x)|^2\chi(x)\, dx=\oint_E|K(x)|^2\chi(x)\, dx.
$$
Therefore
\begin{equation}\label{scprod0}
    \oint_E \overline{K(x)} F(x)\chi(x)\, dx=0.
\end{equation}
We write
$$
K(z)=\frac{k(z)}{k(\infty)}+\tilde K(z),
$$
where, evidently, $\tilde K(\infty)=0$. Putting $F=\tilde K$ in \eqref{scprod0} we get $\|\tilde K\|_\chi=0$.
\end{proof}

Now we use a description \eqref{kade} of reproducing kernels in $H^2$--spaces in Widom domains, for details see \cite{SYU}.

\begin{theorem}\label{threpr} Let $B(z,x_j)$ be the complex Green function, that is, a multi--valued analytic function in $\Omega$ such that $-\log|B(z,x_j)|=G(z,x_j)$. Then
 there exists an  extremal functions of the Problem \ref{pr1}, which  is of  the form
\begin{equation}\label{revil}
    H(z)=-\frac{1}{\sqrt{(z-a_0)(z-b_0)}}\prod_{j\ge 1}\frac{z-x_j}{\sqrt{(z-a_j)(z-b_j)}}\frac{I(\infty)}{I(z)},
\end{equation}
where $I(z)=\prod_j B_{x_j}(z)$  is single--valued in $\Omega$.
\end{theorem}

\begin{proof} Similar to Proposition \ref{pr2.3}, $E^2_\chi$ corresponds   to the $H^2$-space with a suitable character, that is, we can relate the scalar product and the reproducing kernel $K$ in $E^2_\chi$ with the scalar product and the reproducing kernel $K^D$ with respect to the harmonic measure:
\begin{equation*}
    \int_E |K(x)|^2\chi(x)dx=C\int_E|K^D(x)|^2\prod_{j\ge 1}\frac{x-c_j}{\sqrt{(x-a_j)(x-b_j)}}\frac{dx}{\sqrt{(x-a_0)(b_0-x)}}
\end{equation*}
with a certain $D\in D(E)$ and $C>0$. It implies
\begin{equation}\label{zh}
    K(z)=\prod_{j\ge 1}\sqrt{\frac{z-x_j}{\sqrt{(z-a_j)(z-b_j)}}\left(\frac{z-b_j}{z-a_j}\right)^{\tilde \delta_j/2}
    \frac {B_{x_j}(\infty)}{B_{x_j}(z)}}\prod_{j\ge 1} \left(\frac {B_{x_j}(z)}{B_{x_j}(\infty)}\right)^{\frac{1+\epsilon_j}2}.
\end{equation}
Recall that $K(\infty)=1$.
Since $K(z)$ has no zeros in $\Omega$ we have $\epsilon_j=-1$. We substitute \eqref{zh} in  \eqref{exfunc2} to get \eqref{revil}. Since $K(z)$ is singlevalued in $\Omega$ the product $I(z)$ is also singlevalued.

 \end{proof}

\section{Singular components of the reflection\-less \\ measures and DCT}\label{sec5}

\begin{proposition}\label{pr2}
Assume that one of the $R$ functions \eqref{defchi} contains a singular component in its integral representation \eqref{intreprR}.
Then DCT fails.
\end{proposition}

\begin{proof}
Since
 $$
R(z)=\int\frac{d\sigma_{s}(x)}{x-z}+\frac 1 \pi\int_E\frac{|R(x)|dx}{x-z}
$$
we have
$$
\frac 1 \pi\int_E{|R(x)|dx}=1-\sigma_s(E).
$$
Thus $M$ in \eqref{inf} is less than 1 if $\sigma_s(E)>0$.
\end{proof}

The simplest example of a Widom domain where DCT fails was constructed in this way \cite{HA}, see also \cite{NVYu}: assume that
{$E$ is a system of intervals, which accumulate to $b_0$ only.} In this case
there exists a reflectionless measure with a nontrivial masspoint at $b_0$
if and only if
\begin{equation}\label{condm17}
    \int_E\frac{dx}{x-b_0}<\infty.
\end{equation}
So the main point of the example was to demonstrate that \eqref{condm17} does not contradict  the Widom condition.

Note that a much more advanced example   of a reflectionless measure with the \textit{singular continuous component} was constructed in \cite{NVYu}.
Conditions that ensure that all reflectionless measures have no singular component were studied in
\cite{PR}.

Now we show that already   in this simplest case  \eqref{condm17} we have
\begin{equation}\label{smaller}
    M< \inf_{R(z,\{x_j\})}\int_E{|R(x)|dx}.
\end{equation}
Moreover, in the next section we demonstrate that Proposition \ref{pr2} cannot be inverted: the absence of a singular component for all reflectionless  measures does not guarantee DCT for a Widom domain.

To prove \eqref{smaller} we define
\begin{equation}\label{defP}
    \Pi(z)=-\sqrt{\prod_{j\ge 0}\frac{z-b_j}{z-a_j}}=(z-b_0)R(z,\{b_j\}).
\end{equation}
Note that \eqref{condm17} implies that the following limit exists
$$
\lambda_*:=-\lim_{x\uparrow b_0}\Pi(x)=e^{-\frac 1 2\int_E\frac{dx}{x-b_0}}>0
$$
and represents the biggest possible value of the mass for reflectionless measures, i.e.:
\begin{equation}\label{exarextr}
    \inf_{R}\frac 1 \pi\int_E{|R(x)|dx}=1-\lambda_*.
\end{equation}

For $\lambda\in (0,\lambda_*)$ we define
\begin{equation}\label{defRI}
    R_\lambda(z)=\frac 1{1-\lambda^2}\frac{\Pi^2(z)-\lambda^2}{(z-b_0)\Pi(z)},\quad
    I_\lambda(z)=\frac{\Pi(z)+\lambda}{\Pi(z)-\lambda}.
\end{equation}
Here $R_\lambda(z)$ is of the form \eqref{defchi} and $I_\lambda(z)$ is the Blaschke product of the form given in Theorem \ref{threpr}.
For the first function we have
$$
R_\lambda(z)=\frac{\sigma_{0,\lambda}}{b_0-z}+\frac 1\pi\int_E\frac{|R_\lambda(x)|dx}{x-z}.
$$
Moreover
$$
\sigma_{0,\lambda}=\lim_{x\uparrow b_0}(b_0-x)R_\lambda(x)=\frac{\lambda_*^2-\lambda^2}{\lambda_*(1-\lambda^2)}
$$
and therefore
$$
\frac 1\pi\int_E{|R_\lambda(x)|dx}=1-\frac{\lambda_*^2-\lambda^2}{\lambda_*(1-\lambda^2)}=\frac{1-\lambda_*}{\lambda_*}\frac{\lambda_*+\lambda^2}{1-\lambda^2}.
$$

Thus for
$$
H_\lambda(z):=R_\lambda(z)\frac{I_\lambda(\infty)}{I_\lambda(z)}, \quad I_\lambda(\infty)=\frac{1-\lambda}{1+\lambda},
$$
we get
$$
\frac 1\pi\int_E{|H_\lambda(x)|dx}=
\frac {|I_\lambda(\infty)|}\pi\int_E{|R_\lambda(x)|dx}=
\frac{1-\lambda_*}{\lambda_*}\frac{\lambda_*+\lambda^2}{(1+\lambda)^2}.
$$
The last function decreases with $\lambda$, so the smallest value corresponds to $\lambda=\lambda_*$.
On the other hand the extremum \eqref{exarextr} corresponds to $\lambda=0$.
Thus
\begin{equation}\label{exa0}
M_1=\frac 1\pi\int_E{|H_{\lambda_*}(x)|dx}=\frac{1-\lambda_*}{1+\lambda_*}
\end{equation}
and
\begin{equation}\label{exa}
    \inf_{R}\frac 1\pi\int_E{|R(x)|dx}=\frac 1\pi\int_E{|H_0(x)|dx}=1-\lambda_*>
    M_1.
\end{equation}

\begin{remark}
Let
\begin{equation}\label{smcond}
\int_E\frac 1{(x-b_0)^n}dx<\infty.
\end{equation}
In this case the extremum is less than \eqref{exa0}. This smaller value $M_n$ can be expressed by means of a suitable finite (depending on $n$) moment problem.

\end{remark}

\section{No DCT: an infinite dimensional defect space related to a single singular point}

We start with the trace of the following matrix function
\begin{equation}\label{ma1}
    \fA(z)=\begin{bmatrix}
    \cos t\sqrt{z}&\frac{\sin t\sqrt{z}}{\sqrt{z}}\\
    -{\sqrt{z}}\sin t\sqrt{z}&\cos t\sqrt{z}
    \end{bmatrix}\begin{bmatrix}
    \cos {z}&{\sin z}\\
    -\sin {z}&\cos {z}
    \end{bmatrix}, \ t>0.
\end{equation}
It is easy to check that $\det\fA(z)=1$ and
\begin{equation}\label{ma2}
    \frac{\fA^*(z)J\fA(z)-J}{z-\bar z}\ge 0,\quad J=\begin{bmatrix}
    0&-1\\
  1&0
    \end{bmatrix}.
\end{equation}
For this reason $\Delta(z)=\frac 1 2 \tr\fA(z)$ is an entire function with the real $\pm 1$ points. Moreover the function
\begin{equation*}
    \lambda(z)=\frac{\Delta(z)+\sqrt{\Delta(z)^2-1}}{2}
\end{equation*}
is well defined in the upper half-plane, see e.g. \cite{YU}. This is the eigenvalue of $\fA(z)$ with the characteristic property $|\lambda(z)|>1$. So, if we define the domain $\bbC\setminus E_1$,
$E_1=\{x\in R: |\Delta(x)|\le 1\}$, then $\log\lambda(z)$ represents the complex Martin function of this domain with the singular point at infinity. In \cite[VIII]{KU} $\log|\lambda(z)|$ is called the Phragm\'{e}n--Lindel\"of function of the domain.

To simplify further consideration we define  $E=E_1\cup[0,\infty)$.

\begin{lemma}\label{lemma 6.1}
For the given $E$ the domain $\bbC\setminus E$ is of Widom type.
\end{lemma}

\begin{proof}
First of all we use Theorem \cite[p. 407]{KU}. Since $\log|\lambda(iy)|=|y|+o(|y|)$, $y\to\pm\infty$, according to this theorem
\begin{equation*}\label{ku1}
    \int_{-\infty}^\infty  G_1(x,i) dx<\infty,
\end{equation*}
where $G_1(z,i)$ is the Green function of the domain $\bbC\setminus E_1$. Thus  for the extended  boundary $E$ we have also
 \begin{equation}\label{ku2}
    \int_{-\infty}^\infty  G(x,i)dx<\infty,
\end{equation}
where $G(z,i)$ is the Green function of the domain $\bbC\setminus E$.

Using the explicit formula for
\begin{equation}\label{defdelta}
    \Delta(z)=\cos t\sqrt{z}\cos {z}-\frac{z+1} 2 \frac{\sin t\sqrt{z}}{\sqrt{z}}\sin z
\end{equation}
we conclude that on the negative half-axis $E$ consists of a system of intervals $[b_{k+1},a_{k}]$ close to the points $-k\pi$ (the leading term in asymptotics) of the length $a_{k}-b_{k+1}\sim\frac{e^{-t\sqrt{k\pi}}}{\sqrt{k\pi}}$.

Comparing in the standard way \eqref{ku2} with the common area of triangles built on the intervals $(a_k,b_k)$
with the vertex $(c_k,G(c_k,i))$, where $c_k$ is the critical point of $G(z,i)+G(z,-i)$, we get
\begin{equation*}
    \sum G(c_k,i)\frac{b_k-a_k} 2<\infty.
\end{equation*}
Since $b_k-a_k\ge\delta>0$ (in fact, $b_k-a_k \to \pi$) we obtain the Widom condition
$\sum G(c_k,i)<\infty$.
\end{proof}

\begin{proposition}\label{pr6.2} For $\Delta(x)$ given by \eqref{defdelta},
let the system of intervals $(a_k,b_k)$ be defined by the condition $|\Delta(x)|>1$ on the negative half--axis. Then the DCT does not hold in
the Widom domain
$\Omega=(\bbC\setminus \bbR)\cup_k(a_k,b_k)$.
\end{proposition}

\begin{proof} Let us chose a normalization point in a gap, say,  the critical point $c_1\in (a_1,b_1)$. We claim that the function
$F(z)=\cos\tau\sqrt{c_1}-\cos\tau\sqrt{z}$, $0<\tau\le t$,
belongs to the class $E^1_0$ (with respect to $c_1$), that is, it is of Smirnov class in
$\Omega$ and
\begin{equation}\label{nondct}
    \int_E\frac{|F(x)|}{|x-c_1|^2}dx<\infty.
\end{equation}
If so, by the Cauchy Theorem,
\begin{equation*}
    \frac 1{2\pi i}\oint_{\partial\Omega} \frac{F(z)}{(z-c_1)^2}dz=\tau\frac{\sin\tau\sqrt{c_1}}{2\sqrt{c_1}}\not=0,
\end{equation*}
but, in fact, this integral vanishes just due to the symmetry  $F(x+i0)=F(x-i0)$.

We note that $e^{i\tau\sqrt{z}}$ is in absolute value less than one in $\Omega$ and does not vanish. Since the Phragm\'{e}n--Lindel\"of function of the domain behaves as $\Im z$ at infinity this function also does not have the singular inner factor (with the only possible singular point at infinity). Thus it is an outer function in the domain. Therefore $\cos \tau\sqrt{z}$ is of Smirnov class in $\Omega$.

The size of the intervals $[b_{k+1},a_k]$ guaranties that the integral  \eqref{nondct} converges. The proposition is proved.
\end{proof}

\begin{remark}
Note that in the current example $\int_E\frac{dx}{1+|x|}=\infty$, compere \eqref{condm17}. Since infinity is the only possible support for a mass-point,  there is no reflectionless measure with  a singular component, see Sect. \ref{sec5}. Evidently, this set is weakly homogeneous, so we can also refer to the general result \cite{PR}. Thus
 $\inf_{R}\frac 1 \pi\int_E{|R(x)|dx}=1$, but DCT fails and $M=M(E)<1$ in Problem \ref{pr1}.
\end{remark}

\begin{remark}
Let us mention that the example in \cite{HA}, which was discussed in Sect. \ref{sec5}, corresponds to the case when the defect space $H^2(\alpha)\ominus \check H^2(\alpha)$ (for a certain $\alpha$) has dimension  one (in particular, it is non--trivial). The example with a singular measure \cite{NVYu} corresponds to an infinite defect space, but with infinitely many "singular" points in the domain. Our example corresponds to \textit{an infinite-dimensional defect space related to a single singular point} in the domain.
This remark explains the main idea of the current constuction. The form of the product \eqref{ma1} is dictated by a  simple reason: each defect space  generates a factor with a smaller grow nearby the singular point comparably with the grow of the Martin function of the corresponding domain (in the given case $O(\sqrt{z})$ and $O(z)$ respectively).
We will discuss such relations in details in a forthcoming paper.
\end{remark}

\section{Widom domain with DCT and \\ non-homogeneous boundary}\label{secbndx}
First of all we note that  Theorem \ref{threpr} was stated for a bounded set $E$. We will consider a domain with a unique accumulation point, say $b_0$, for the ends of intervals $a_k$'s and $b_k$'s. It is convenient to send this point to infinity by the change of the variable $z_1=\frac{1}{b_0-z}$. In this case a function with the only possible singularity at $b_0$ becomes a standard entire function.

Let $E$ be an unbounded closed set.
Note that an unbounded closed set $E$  is homogeneous if \eqref{hocond} holds for all $x\in E$ and \textit{for all} $\delta>0$.
 In \cite{SYU} it was shown that homogeneity of the boundary of the Denjoy domain implies DCT. In this section we show that homogeneity is not a necessary condition for DCT.

\begin{lemma}\label{lemma 7.1} Let $E$ be a system of intervals on the negative half-axis which accumulate to infinity only.
Assume that $\Omega=\bbC\setminus E$ is of Widom type and DCT fails. Then there exists an entire function
$F$ of Smirnov class in $\Omega$ such that
\begin{equation}\label{efsod}
    \int_E |F(x)|\frac{dx}{|x|}<\infty.
\end{equation}

\end{lemma}

\begin{proof} If $\int_E \frac{dx}{|x|}<\infty$ then $F(x)=1$. If not, then there is no reflectionless measure  with
a singular component.  Thus, the factor $I(z)$ in Theorem \ref{threpr}  is not a constant. We define
\begin{equation}\label{deffct}
    F(z)=\frac{I(0)}{\sqrt{1-z/a_0}}\prod_{j\ge 1}\frac{1-z/x_j}{\sqrt{(1-z/a_j)(1-z/b_j)}}\left\{\frac{1}{I(z)}-I(z)\right\}.
\end{equation}
It is of Smirnov class in $\Omega$, moreover it is real valued on the whole real axis. Since
$ \frac{1}{I(x\pm i0)}-I(x\pm i0)=2i\Im \frac{1}{I(x\pm i0)}$, the integral \eqref{efsod} converges and all possible singularities on $E$ are removable. Thus $F(z)$ is an entire function.
\end{proof}

\begin{remark}\label{rem6.2}
 Recall that $|I(z)|<1$ in the domain. Thus $F(z)$ has only real zeros, moreover only in  the set $E$.
 Note also that as soon as $F(z)$ is not a constant there exists $F_1(z)$ such that
 $\int_E |F_1(x)|dx<\infty$. For instance $F_1(z)=\frac{F(z)}{z-x_0}$, where $x_0$ is a zero of $F(z)$.
 \end{remark}

We denote by $\cM_E(z)$ the Martin function for $\bbC\setminus E$ with singularity at infinity.
The function $\cM_E(z)$ has at most order 1/2, that is, $\cM_E(z)=O(\sqrt{|z|})$.
 We say that a set $E$ is an Akhiezer-Levin set \cite{AkhL, BoS2} if
 \begin{equation*}
    \limsup_{|z|\to\infty} \frac{\cM_E(z)}{\sqrt{|z|}}>0.
 \end{equation*}
 It is worth mentioning that in this case %the limit
 %\begin{equation*}
 $\lim_{x\to+\infty} \frac{\cM_E(x)}{\sqrt{x}}$
 %\end{equation*}
 exists.

  Using the change of variable $z=-z_1^2$ we can work with even functions in the upper half-plane, related, correspondingly, to  symmetric subsets of the real axis. In particular, in this case we can refer directly to
  \cite[Sect. VIII]{KU}.

 \begin{theorem}\label{thdctyes} Let $E=\bbR\setminus \cup_{k=-\infty}^\infty(a_k,b_k)$, $a_{- k}=- b_k$,
 consist of uniformly separated intervals, say,
 \begin{equation}\label{unsep}
    b_k-a_k\ge 1,\quad\text{for all}\ k.
 \end{equation}
Assume  that the following "weighted" Widom condition holds true
\begin{equation}\label{modwid}
    \sum_{c_k\in (a_k,b_k):G'(c_k,i)=0}G(c_k,i)(b_k-a_k)<\infty.
\end{equation}
For an arbitrary $\delta\in (0,1/2)$ define
\begin{equation}\label{defedel}
    E_\delta=\bbR\setminus \cup_{k=-\infty}^\infty(a_k+\delta,b_k-\delta).
\end{equation}
If
\begin{equation}\label{intcond}
    \int_{E_\delta}\frac{dx}{1+|x|}=\infty,
\end{equation}
then the domain $\Omega_\delta=\bbC\setminus E_\delta$ is of Widom type with DCT.
\end{theorem}

\begin{proof}
First of all \eqref{unsep} and \eqref{modwid} imply that $\bbC\setminus E$ is of Widom class.
Further, \eqref{modwid} implies $\int_{\bbR} G(t,z) dt<\infty$. By the Koosis' criterion \cite{KU}  $E$
is a symmetric Akhi\-ezer-Levin set, that is,
  \begin{equation*}
 \lim_{y\to+\infty} \frac{\cM_E(iy)}{{y}}>0.
 \end{equation*}
Evidently the extended set $E_\delta$ \eqref{defedel} belongs to the Akhiezer-Levin class and $\Omega_\delta$ is of Widom class.

Assume that DCT fails. By Lemma \ref{lemma 7.1} there exists a  non trivial even entire   function $\tilde F(z)$ of Smirnov class in $\Omega_\delta$ such that
\begin{equation*}
    \int_{E_\delta} |\tilde F(x)|\frac{dx}{1+|x|}<\infty.
\end{equation*}
Due to \eqref{intcond} $\tilde F(z)$ is not a constant and we can find a \textit{non trivial}  entire   function $F(z)$, see  Remark \ref{rem6.2},
 such that
\begin{equation}\label{boint}
    \int_{E_\delta} |F(x)|{dx}\le 1.
\end{equation}
Since $F(z)$ is of Smirnov class in $\Omega_\delta$ we have
\begin{equation*}
    \lim_{|z|\to+\infty}\frac{\log|F(z)|}{\cM_{E_\delta}(z)}=0,\quad z=x+iy,\ y\ge c|x|,\ c>0.
\end{equation*}
Moreover, since $F(z)$,  in  particular, is in the Cartwright class and all its zeros are real, one can show, see e.g. \cite[p. 58]{AKH}, that for all $\epsilon>0$ the following a priori estimate holds true $|F(z)|\le A(\epsilon)e^{\epsilon|z|}$ in the whole complex plane.

Now for $\delta'<\delta$ we define the entire function
\begin{equation*}
    H(z)=H(z,\delta')=\int^{\delta'}_{-\delta'} F(z+t) dt.
\end{equation*}
By \eqref{boint} $H(z)$ is uniformly bounded  on $E$, $|H(x)|\le 1$, $x\in E$, and also possesses the a priori estimate $|H(z)|\le A(\epsilon)e^{\epsilon|z|}$. We apply the Phragm\'{e}n--Lindel\"of principle
to $H(z)$ in the domain $\Omega=\bbC\setminus E$, see \cite[p. 406]{KU}, to get that in fact
$|H(z)|\le e^{\epsilon\cM_E(z)}$. Since $\epsilon$ is arbitrary small, $H(z)=H(z,\delta')$ is bounded in the whole complex plane, therefore it is a constant. Since this holds for an arbitrary positive $\delta'<\delta$,  $F(z)$ is a constant, and due to \eqref{boint} and \eqref{intcond} $F(z)=0$.
The contradiction with the claim that $F(z)$ is a non-trivial function shows that DCT holds in $\Omega_\delta$.

\end{proof}

\begin{example}[Benedicks' set] Let $p>1$ and put
\begin{equation*}
    E_\delta=\cup_{n=1}^\infty\{[-n^p-\delta,-n^p+\delta]\cup[n^p-\delta,n^p+\delta]\},
\end{equation*}
$\delta>0$ being taken small enough so that the intervals figuring in the union do not intersect.
This set is not homogeneous and it is of Akhiezer-Levin class
 due to the Koosis' criterion.   In fact, due to Benedicks \cite[p. 439]{KU}
$$
G(x,i)\le C(\delta)\frac{\log|x|}{|x|^{(p+1)/p}},
$$
i.e., it is integrable on the real axis and \eqref{modwid} is satisfied.

The set has a finite logarithmic length \eqref{intcond}. So we modify it by adding each second interval formed by the geometric progression:
\begin{equation*}\label{benmod}
    E_{\delta,q}=E_\delta\cup_{n=1}^\infty\{[-q^{2n},-q^{2n-1}]\cup[q^{2n-1},q^{2n}]\}, \ q>1.
\end{equation*}
 The resulting set is still  not homogeneous, but satisfies all requirements of Theorem \ref{thdctyes} (we can add to the set all possible gaps of the standard length less than one, in case such open intervals are present in $\bbR\setminus E_{\delta,q}$).
Thus the DCT holds in $\Omega:= \bbC\setminus  E_{\delta,q}$.
\end{example}

\begin{corollary} For $p>1$  let
$$
\tilde E=\{0\}\cup\{x:1/x\in E_{\delta,q}\},
 \ \delta\in(0,1), \ q>1.
$$
Then every reflectionless Jacobi matrix $J\in J(\tilde E)$ is almost periodic.
\end{corollary}

%section{attachment}
%
%Let $z=-x$, then we have
%$$
%\Delta(-x)=\cosh t\sqrt{x}\cos {x}-\frac{x-1} 2 \frac{\sinh t\sqrt{x}}{\sqrt{x}}\sin x=\pm 1
%$$
%Define
%\begin{equation*}
%    A(x)=\sqrt{\left(\frac{x-1} 2 \frac{\sinh t\sqrt{x}}{\sqrt{x}}\right)^2+(\cosh t\sqrt{x})^2}\sim\sqrt{x}e^{t\sqrt{x}}
%    \end{equation*}
%and
%\begin{equation*}
%    \sin \phi(x)=\frac{\cosh t\sqrt{x}}{A(x)}\sim \frac{1}{\sqrt{x}}
%\end{equation*}
%In this case we have
%\begin{equation*}
%    \sin(x+\phi(x))=\pm\frac 1{A(x)}.
%\end{equation*}
%That is
%\begin{equation*}
%    a_{k+1}-\phi(a_{k+1})=2\pi k+\frac 1{A(a_{k+1})}
%\end{equation*}
%and
%\begin{equation*}
%    b_{k}-\phi(b_{k})=2\pi k-\frac 1{A(a_{k+1})}
%\end{equation*}
%Thus
%\begin{equation*}
%    a_{k+1}-b_{k}-(\phi(a_{k+1})-\phi(b_{k}))=\frac 1{A(a_{k+1})}+\frac 1{A(b_{k})}
%\end{equation*}
%Due to the given asymptotics
%\begin{equation*}
%    (a_{k+1}-b_{k})(1+\frac{1}{\sqrt(a_{k+1})\sqrt(b_{k})(\sqrt(a_{k+1})+\sqrt(b_{k}))}=\frac 1{A(a_{k+1})}+\frac 1{A(b_{k})}
%\end{equation*}
%Finally
%$a_{k+1}-b_{k}\sim\frac{e^{-t\sqrt{k\pi}}}{\sqrt{k\pi}}$.

\bibliographystyle{amsplain}

%\newpage

\bigskip
 Abteilung f\"ur Dynamische Systeme und Approximationstheorie,

 Johannes Kepler University Linz,

A-4040 Linz, Austria

\smallskip
\textit{E-mail address:}

Petro.Yudytskiy@jku.at

\end{document}